\documentclass{article}
\usepackage[utf8]{inputenc}
\usepackage{blindtext}
\usepackage[utf8]{inputenc}
\usepackage{graphicx}  
\usepackage{listings}
\usepackage{tikz}
\usetikzlibrary{arrows}
\usepackage{amsmath}
\usepackage{xfrac}
\usepackage{xspace}
\usepackage{stmaryrd}
\usepackage{mathrsfs}
\usepackage{amsmath, amssymb, amsthm, amsfonts, mathtools,amscd}
\usepackage{graphicx, setspace, subfig}
\usepackage{frontespizio}

\newtheorem{definizione}{Definition}[section]
\newtheorem{teorema}{Theorem}[section]
\newtheorem{proposizione}{Proposition}[section]
\newtheorem{corollario}{Corollary}[section]

\newtheorem{osservazione}{Remark}[section]

\newenvironment{dimostrazione}{\medskip\noindent{\bf Proof}\enspace} {\hfill\newline\smallskip}
\title{On the continuity of Weil-Petersson volumes of the moduli space weighted points on the projective line}
\author{Salvatore Tambasco  \footnote{salvatore.tambasco01@universitadipavia.it}}
\date{}

\begin{document}

\maketitle

\begin{abstract}
   \noindent 
   In this work we show that the Weil-Petersson volume (which coincides with the CM degree) in the case of weighted points in the projective line is continuous when approaching the Calabi-Yau geometry from the Fano geometry. More specifically, the CM volume computed via localization converges to the geometric volume, computed by McMullen with different techniques, when the sum of the weights approaches the Calabi-Yau geometry.
\end{abstract} 

\section{Introduction}

When considering a moduli space of weighted pairs it becomes interesting to study how such moduli space changes towards the \emph{geography} of the weights. In this direction, we study the behaviour of the \emph{volume function} of the moduli space of points in $\mathbb{P}^1$. We denote the latter by $M_d,$ where $d$ denotes the weight vector $(d_1, ..., d_n)$ of a configuration of $n$-points in $\mathbb{P}^1.$ We assume, henceforth, that the components of such vector are rational numbers in the open interval $(0,1).$ A point of $M_d$ is a pair $(\mathbb{P}^1, \sum_{i=1}^n d_i p_i ),$ where $p_i$ are  points in general position in $\mathbb{P}^1,$ $\forall i \{1,2, ...,n\}.$  The sum of the weights determine the nature of the pair. Namely, given the pair $(\mathbb{P}^1 , \sum_{i=1}^n d_i p_i),$ a choice of the rational weights correspond to three distinguished \emph{geometries} 

\begin{displaymath} 
 \begin{cases} 
 \sum_{i=1}^n d_i < 2 & \ \text{log Fano,} \\
 \sum_{i=1}^n d_i =2 & \ \text{log Calabi-Yau,} \\
 \sum_{i=1}^n d_i > 2 & \ \text{log General type.}
 \end{cases}
 \end{displaymath}
 
We want to study how the Weil-Petersson volume (induced by \emph{canonical constant curvature metrics}) of the moduli space of $n$ points in the projective line varies in the above geometries, henceforth denoted by  $M_d.$ This study is motivated by the relations among the first Chern class of the CM line bundle and the Weil-Petersson metric. The aim of this work is to show that the CM volume of $M_d$ \emph{converges} to the geometric volume computed by McMullen in \cite{MM00} (Theorem 8.1), when the sum of the weights approaches the log Calabi-Yau geometry. To do so, we take the below explicit description of $M_d$ as GIT quotient

 \begin{equation} \label{smallweight} M_d =  (\mathbb{P}^1)^{n} //_{\mathcal{L}_a} \mathrm{SL}(2),  \end{equation}

where $d = (d_1, ..., d_n)$ is a weight vector. This compact space parametrizes configurations of points in $\mathbb{P}^1$ which carry canonical metrics. Indeed, when  $\sum_{i=1}^n d_i = 2,$ it parametrizes flat metrics (with cone angle) thanks to the work of \cite{Th98}, \cite{DeliMosto}. When $\sum_{i=1}^n d_i < 2,$  by works of Li in \cite{LX14} and Fujita in \cite{Fuji19b}, the latter moduli space is actually a Fano K-moduli space, parametrizing K-polystable Fano pairs (positive constant curvature metric). Finally,  when $\sum_{i=1}^n d_i = 2 + \epsilon,$ and generic then it is a Hasset/KSBA type moduli space (\cite[Theorem 5.5.2]{Ale15}, \cite[Theorem 1.5]{Ale08}) of constant negative curvature metric.

Volumes of $M_d$ have been computed in many works, via localisation by \cite{Man08}, \cite{Ta01}, \cite{Khoi05}, \cite{Mar00}, using the isomorphism with the moduli space of polygons. The correspondence between the Weil-Petersson volume and the degree of the CM line bundle it is shown, for the Fano case, in \cite{LWX15} in the absolute case, and the latter result is extended in the log case in \cite{Tam20}.\\

In this work we will prove the following:

\begin{teorema} \label{importante} The Fano CM volumes of $M_d$ with $d_i \in (0,1) \cap \mathbb{Q}$, for all $i \in \{1,2,....,n\}$ converge to the volume computed by McMullen in \cite{MM00} (Theorem 8.1) when the sum of the weights approaches 2 from below.
\end{teorema}

\begin{osservazione} In \cite{MM00} McMullen computed the volume of the complex hyperbolic metric introduced by Thurston in \cite{Th98} on $M_d$ for $\sum_{i=1}^n d_i = 2$.   
\end{osservazione}

The above proves the continuity of the \emph{volume function} when passing from the Fano geometry to the Calabi-Yau geometry. The continuity result holds also when the sum of the weights is $2 + \epsilon,$ due to \cite[Theorem 5.5.2]{Ale15}, \cite[Theorem 1.5]{Ale08}. When the sum of the weights exceed significantly two, then we do not have the above explicit description of the moduli space of weighted points in $\mathbb{P}^1.$ Namely, the natural K-moduli compactification is not simply given by the above GIT quotient. However, in this work, as an example, we fully examine the simple case of four points in $\mathbb{P}^1$ giving an explanation of the continuity of the volume function across the above mentioned geometries through the picture of $K$-stability.

We organised this work as follows

\begin{itemize}
    \item In section 2 we recall the notion of CM line bundle for pairs and we relate it to the moduli of log Fano hyperplane arrangements.
    \item In section 3 we fully examine the case of four points on the projective line, by calculating the volume of the moduli space of four points in $\mathbb{P}^1.$ We prove Theorem \ref{importante} in this case. Later we discuss the continuity of the volume function across the above mentioned three geometries.
    \item In section 4 we give a proof of Theorem \ref{importante} by using a continuity argument arising from the continuity of the construction of the Weil-Petersson metrics.
    \item In the final remark, as a \emph{miscellanea} we show that the CM volume can be expressed as the sum of Donaldson-Futaki invatiant of some \emph{special} test configurations.
\end{itemize}

\section*{Acknowledgments}

I would like to express my most sincere gratitude towards my supervisors Prof. Cristiano Spotti and Prof. Alessandro Ghigi for their patient guidance and advice they have provided me. Last, but not least I wish to extend my thanks to the Mathematics Department of Aarhus University for hosting me for two years, to the AUFF starting grant 24285, and to the Villum Fonden 0019098 for partially funding my research.

\section{The log CM line bundle}

\subsection{The CM line bundle for pairs}

Fix $f : \mathcal{X} \rightarrow B$ to be a proper flat morphism of scheme of finite type over $\mathbb{C}.$ Let $\mathcal{L}$ be a relatively ample on $\mathcal{X}$ and assume that $f$ has relative dimension $n \geq 1,$ namely $\forall b \in B, \ (\mathcal{X}_b, \mathcal{L}_b)$  has constant dimension $n.$ Below, in order to follow the conventions used in the \emph{area}, multiplicative and additive notation for divisors and line bundles are often interchanged. 

\begin{definizione} \label{ad} \cite{ADL19} Let $\mathcal{D}_i, \forall i=1,2, ..,k $ be a closed subscheme of $\mathcal{X}$ such that $f|_{\mathcal{D}_i} : \mathcal{D}_i \rightarrow B$ has relative dimension $n-1,$ and $f|_{\mathcal{D}_i}$ is proper and flat. Let $d_i \in  [0,1] \cap \mathbb{Q}$ we define the log CM $\mathbb{Q}$-line bundle of the data $(f, \mathcal{D}:= \sum_{i=1}^k d_i \mathcal{D}_i), \mathcal{L})$ to be 

$$ \lambda_{CM, \mathcal{D}} := \lambda_{CM} - \frac{n \mathcal{L}_b^{n-1} \cdot \mathcal{D}_b}{(\mathcal{L}_b^n)} \lambda_{\mathrm{CH}} + (n+1) \lambda_{\mathrm{CH}, \mathcal{D}}, $$  
Where $\lambda_{CM}$ is the CM line bundle defined in Chapter two and $\lambda_{CH}$ is the Chow line bundle, defined as the leading order term of the Hilbert Mumford expansion, namely $\lambda_{CH}= \lambda_{n+1},$ and $\lambda_{CH, \mathcal{D}} = \bigotimes_{i=1}^k \lambda_{CH}^{d_i}.$
\end{definizione}

We outline some easy and well known properties of the log CM line bundle, for a full description of the below facts see \cite{ADL19}.

\begin{proposizione} \label{pr} Let $f : (\mathcal{X}, \mathcal{D}) \rightarrow B$ be a $\mathbb{Q}-$Gorenstein flat family of $n$ dimensional pairs over a normal proper variety $B.$ Then for any $\mathcal{L}$ relatively $f-$ample line bundle we have

\begin{equation} \label{cmlog} c_1(\lambda_{CM, \mathcal{D}}) = n \frac{ (-K_{\mathcal{X}_b} - \mathcal{D}_b) \cdot \mathcal{L}_b^{n-1}}{(\mathcal{L}_b^n)} f_{*} c_1 ( \mathcal{L})^{n+1} - (n+1)f_{*}((-K_{\mathcal{X}/B} - \mathcal{D}) \cdot c_1(\mathcal{L})^n). 
\end{equation}
\end{proposizione}

\begin{osservazione} In \cite{SGM18}, there is another definition of the log CM line bundle for the case of one divisor with weight $1-\beta$. And in the same work (Theorem 2.7 first assertion) there is a calculation of the first Chern class of the defined log CM line bundle. We point out that the definition we gave in this section it is exactly the same for the case of one divisor. Indeed, for the first assertion of Theorem 2.7 we just require to let

$$ \mu(\mathcal{L}_b) := - \frac{-K_{\mathcal{X}_b} \cdot \mathcal{L}^{n-1}}{(\mathcal{L})^n}, \  \ \mu(\mathcal{L}_b, \mathcal{D}_b) := \frac{\mathcal{D}_b \cdot \mathcal{L}_b^{n-1}}{(\mathcal{L}^n)}$$ 
and of course $\mathcal{D} = (1- \beta)D.$ 

The definition of the log CM line bundle given in the same work \cite{SGM18} definition 2.2 is also a functorial definition, namely 

$$ \Lambda_{CM, \mathcal{D}} := \lambda_{n+1}^{n(n+1) + \frac{2a_1 - (1- \beta)\tilde{a}_0}{a_0}} \otimes \lambda_n^{-2(n+1)} \otimes \tilde{\lambda}_n^{(1 - \beta)(n+1)} $$
The element $\tilde{\lambda}_n$ refers to the leading order term of the Knudsen-Mumford expansion for $f_{!|\mathcal{D}} \mathcal{L}$, and $\tilde{a}_0$ to the leading order term coefficient of the corresponding Hilbert polynomial, that is by the Hirzebruch-Riemann-Roch expansion

$$ \tilde{a}_0= \frac{1}{n!} \mathcal{L}_b^n \cdot \mathcal{D}_b.$$ 

Unravelling the other terms via the Hirzebruch-Riemann-Roch expansion, we obtain definition \ref{ad}

\end{osservazione}

When computing the intersection number of the log CM line bundle on a $\mathbb{Q}-$Gorenstein flat family on the moduli of pairs, it is important to choose a ample line bundle on such family. In the state of the art of these possible choices, mainly we can make two of such. These choices are presented in the following

\begin{corollario}\label{cpsgm} Let $f: (\mathcal{X}, \mathcal{D}) \rightarrow B$ be a $\mathbb{Q}-$Gorenstein flat family over a normal proper variety $B.$  We have

\begin{enumerate}

\item \cite{CodPat} Given $\mathcal{L}= -K_{\mathcal{X}/B} - \mathcal{D}$, then $$ c_1(\lambda_{CM, \mathcal{D}}) = -f_{*}c_1(-K_{\mathcal{X}/B} - \mathcal{D})^{n+1}.$$
\item \cite{SGM18} Given $\mathcal{L}= -K_{\mathcal{X}/B},$ and $\mathcal{D}_{\mathcal{X}_b} \in | -K_{\mathcal{X}_b} |, \ \forall b \in B$, then $$ c_1 (\lambda_{CM, \mathcal{D}})= f_{*}\left( c_1 (-K_{\mathcal{X}/B})^{n} \cdot (-c_1 (-K_{\mathcal{X}/B}) + \sum_{i=1}^k d_i((n+1) \mathcal{D}_i - n c_1 (-K_{\mathcal{X}/B})) \right). $$ \end{enumerate}
\end{corollario}

Depending on the choice of $\mathcal{L}$ we will get different results in terms of intersection number on a generic curve into the moduli of pairs we want to study. In the below subsection, as an example, we will point out this difference when calculating the intersection number of the log CM line bundle with a generic $\mathbb{P}^1$-\emph{curve} in the moduli space of log Fano hyperplane arrangements.
Finally, we wish to make a remark on the log CM line bundle in the case of families of general type varieties. As a direct consequence of Proposition \ref{pr}, we get the following 

\begin{corollario} \label{loggtcr} Let $f \colon (\mathcal{X}, \mathcal{D}) \rightarrow B$ be a $\mathbb{Q}$-Gorenstein family of log general type varieties of relative dimension n, and with relatively ample line bundle $\mathcal{L}= K_{\mathcal{X}/B} + \mathcal{D}.$ Then,

\begin{equation} \label{logcmgt} 
c_1(\lambda_{\mathrm{CM}}) = f_{*}c_1 (K_{\mathcal{X}/B} + \mathcal{D})^{n+1}
\end{equation}
\end{corollario}

\begin{dimostrazione} By applying directly Proposition \ref{pr} we get

$$ n\frac{(-K_{\mathcal{X}_b} - \mathcal{D}_b) \cdot (K_{\mathcal{X}_b} + \mathcal{D})^{n-1}}{(K_{\mathcal{X}_b} + \mathcal{D}_b)^n} = -n. $$Thus,

\begin{align*} c_1(\lambda_{\mathrm{CM}, \mathcal{D}}) &= -n f_{*}c_1(K_{\mathcal{X}/B} + \mathcal{D})^{n+1} - (n+1)f_{*}(-(K_{\mathcal{X}/B} + \mathcal{D})c_{1}(K_{\mathcal{X}/B} + \mathcal{D})^{n}) \\
&= -n f_{*}c_1( K_{\mathcal{X}/B} + \mathcal{D})^{n+1} + (n+1) f_{*}c_1 (K_{\mathcal{X}/B} + \mathcal{D})^{n+1} \\
&= f_{*}c_1(K_{\mathcal{X}/B} + \mathcal{D})^{n+1}.
\end{align*}
\hfill $\Box$
\end{dimostrazione}

\subsection{The intersection number of the log CM line bundle}

In this section we compute the intersection number of the log CM line bundle with a family of log Fano hyperplane arrangements with base $\mathbb{P}^1.$ Consider the product $\mathcal{X} := \mathbb{P}^n \times \mathbb{P}^1,$ We define a divisor $ \mathcal{D} \subset \mathcal{X},$ as 
$$\mathcal{D} := \sum_{i=1}^{m-1} d_i \mathrm{pr}_2^{*} h + d_m (\mathrm{pr}_1^{*}l + \mathrm{pr}_2^{*}h)$$
Where the $\mathrm{pr}_i$ denotes the projections onto the first and second factor of $\mathcal{X}.$ We can think about $\mathcal{D}$ as $m-1$ fixed hyperplanes of $\mathbb{P}^n$ together with a line $l$ with weight $d_m$free to move along the \emph{diagonal} of $\mathcal{X}.$ We assume that $d_i \in (0,1) \cap \mathbb{Q}, \ \forall i \in \{1,2, .., m\},$ $\sum_{i=1}^{m} d_i < n+1,$ and the chosen line and hyperplanes are in general position. 

\begin{proposizione} \label{prop1} With the above data we have

$$c_1(\lambda_{CM, \mathcal{D}}) = (n+1)d_j \left(n+1 - \sum_{i=1}^m d_i \right)^n, \ \forall j \in \{1,2, ...,m\}. $$
\end{proposizione}

\begin{dimostrazione} We see that the projection $\mathrm{pr}_2 : \mathcal{X} \rightarrow \mathbb{P}^1$ gives a proper and flat family of relative dimension $n.$ Choose $\mathcal{L}= -K_{\mathcal{X}/B} - \mathcal{D}$ and use assertion $1.$ of \ref{cpsgm}, we find that
$$ \mathcal{L}= \left( n+1 - \sum_{i=1}^m d_i \right) \mathrm{pr}_1^{*}h - d_m \mathrm{pr}_2^{*}l $$
Using the binomial expansion, the only term that survives is 

$$c_1(\lambda_{CM,\mathcal{D}}) = -\mathrm{pr}_{2*} c_1 (\mathcal{L})^{n+1} = (n+1)d_m  \left( n+1 - \sum_{i=1}^m d_i \right)^n \mathrm{pr}_{2*} \left( \mathrm{pr}_1^{*} h \cdot \mathrm{pr}_2^{*}l \right). $$ By using the projection formula and the arbitrariness of the choice of the weight on the line, the claim follows. \hfill $\Box$

\end{dimostrazione}

\begin{osservazione} If we would have chosen $\mathcal{L}= - K_{\mathcal{X}/B},$ we can see that the hypothesis $\mathcal{D}_{|\mathcal{X}_b} \in |-K_{\mathcal{X}_b}|$ would not be satisfied. Indeed, it is true if and only if $\sum_{i=1}^m d_i = (n+1),$ i.e. a Calabi-Yau log hyperplane arrangement, but this is in contradiction with the log Fano assumption. However, $-K_{\mathcal{X}/B}$ it is still a relatively ample line bundle on $\mathrm{pr}_2 : \mathcal{X} \rightarrow \mathbb{P}^1.$ We can compute the CM line bundle with this choice, and relaxing the hypothesis for which $\mathcal{D}_{|\mathcal{X}_b} \in |-K_{\mathcal{X}_b}|$  as the next result will show. 
\end{osservazione}

\begin{proposizione} \label{prop2} With the above data and $\mathcal{L}= -K_{\mathcal{X}/B},$ we have

$$ c_1 (\lambda_{CM, \mathcal{D}})=(n+1)^2 d_j. $$\end{proposizione}

\begin{dimostrazione} The first summand is zero, since $\mathcal{X}/\mathbb{P}^1$ is a trivial fibration. From the second summand we get
\begin{multline*}
  c_1 (\lambda_{CM, \mathcal{D}}) = \\
  =-(n+1)\mathrm{pr}_{2*} \left[ \left( n+1 - \sum_{i=1}^m d_i \right) \mathrm{pr}_2^{*}h - d_m \mathrm{pr}_{1}^{*}l \right] \cdot (n+1) \mathrm{pr}_{2}^{*}h= \\
  =(n+1)^2 d_m.
\end{multline*}
Because of the arbitrariness of choices the claim follows $\forall j \in \{1,2, ..., m\}.$ \hfill $\Box$

\end{dimostrazione}

From the choice of the two ample line bundles, see corollary \ref{cpsgm}, on the given flat and proper family we see that they mainly differ by a factor which involves the volume of the fiber, that is 

$$ \left( n+1 - \sum_{i=1}^m d_i \right)^{n}. $$

This difference might become important when calculating the volume of the moduli space of weighted hyperplane arrangements. To see that, we shall calculate the log CM line on to the mentioned moduli space. Recall that the moduli space of weighted log Fano hyperplane arrangement is the GIT quotient

$$ M_d := (\mathbb{P}^n)^m //_{\mathcal{L}_d} \mathrm{SL}_{n+1} $$
where $\mathcal{L}_d = \mathcal{O}(d_1 , ..., d_m),$ is the linearization and $d_i \in (0,1) \cap \mathbb{Q}, \forall i \in \{1,2, ..., m \}.$ It is known that, 

$$\mathrm{Pic} (M_d) = \mathbb{Z}^m $$ 
see for example \cite[Chapter 11, Lemma 11.1]{Dol03}.  Therefore $c_1 (\lambda_{CM, \mathcal{D}}) \in \mathrm{Pic}(M_d) \Rightarrow c_1(\lambda_{CM, \mathcal{D}})= \mathcal{O}(r_1, ..., r_m),$ for some $(r_1, ..., r_m) \in \mathbb{Z}^m.$ In order to compute the $r_j'$s we can use the results of \ref{prop1} and \ref{prop2}. Namely, in both cases we compute the following intersection number

$$ c_1 (\lambda_{CM, \mathcal{D}}) \cdot (\mathcal{X} \rightarrow \mathbb{P}^1) = k, $$
which means, $\forall j \in \{1,2, ..., m\}$ 

$$ \int_{(\mathbb{P}^1)_j} \mathcal{O}(r_1, ..., r_m) = r_j=k .$$
Therefore, from \ref{prop1} we get $r_j = (n+1)d_j \mathrm{Vol}(\mathrm{pr}_{2,t}),$ and from \ref{prop2} we get $r_j = (n+1)^2 d_j.$ Hence, we just proved the following result.

\begin{proposizione} \label{prop3} The log CM line bundle on the moduli space of log Fano hyperplane arrangement is given by

\begin{itemize}

\item  When $\mathcal{L}= -K_{\mathcal{X}/B} - \mathcal{D},$ we have that $$\lambda_{CM, \mathcal{D}} =\left( n+1 - \sum_{i=1}^m d_i \right)^{n} \mathcal{O}((n+1)d_1, ..., (n+1) d_m)$$

\item  When $\mathcal{L}= -K_{\mathcal{X}/B},$ we have $$\lambda_{CM, \mathcal{D}} =(n+1) \mathcal{O}((n+1)d_1, ..., (n+1)d_m).$$\end{itemize}
\end{proposizione}

\begin{osservazione} \label{remarko} In both cases of \ref{prop3} we can see that the log CM line bundle is a multiple of the \emph{weighted prequantum} line bundle on $(\mathbb{P}^n)^m$hence of its K\"{a}hler form, we recall that 

$$\left[ \frac{\omega}{2 \pi}\right] = c_1(-K_{(\mathbb{P}^n)^m})= c_1(\mathcal{O}((n+1)d_1, .., (n+1)d_m))) $$
Where the K\"{a}hler form $\omega$ is given by the weighted sum of the Fubini Study metrics on each $\mathbb{P}^n$, namely $\omega = \sum_{i=1}^m d_i \omega_{FS,i}.$  
\end{osservazione}

\section{A first study: the case of four points on the complex projective line.}

As a first study, we consider the case of four points $\{p_i \}_{i=1}^4$ in $\mathbb{P}^1$ with rational weights $d_i \in (0,1) \cap \mathbb{Q}, i \in \{ 1,2,3,4\}.$ Given the pair $(\mathbb{P}^1 , \sum_{i=1}^4 d_i p_i),$ a choice of the rational weights correspond to three distinguished \emph{geometries}

\begin{displaymath} 
 \begin{cases} 
 \sum_{i=1}^4 d_i < 2 & \ \text{log Fano,} \\
 \sum_{i=1}^4 d_i =2 & \ \text{log Calabi-Yau,} \\
 \sum_{i=1}^4 d_i > 2 & \ \text{log General type.}
 \end{cases}
 \end{displaymath}
 
We want to study how the \emph{Volume} of the moduli space of four points in the projective line varies in the above geometries. We begin with the Calabi-Yau geometry. The moduli space of four points in $\mathbb{P}^1$ can be described as the following GIT quotient

\begin{equation} \label{git4} M_d = (\mathbb{P}^1)^4 //_{\mathcal{L}_d} \mathrm{SL}_2 \mathbb{C}, \end{equation}
with linearization $\mathcal{L}_d =\mathcal{O}(d_1, d_2, d_3, d_4).$ Furthermore, $M_d$ is isomorphic to the moduli space of marked curves of genus $0,$ denoted by $\mathcal{M}_{0,4}$. For clarity, we recall the following result 

\begin{teorema} (\cite{MM00}, Theorem 8.1) \label{mc2} Let $\mathcal{M}_{0,n}$ be the moduli space of $n$ ordered points on the Riemann sphere. Then the complex hyperbolic volume of $\mathcal{M}_{0,n}$ is given by
 
 \begin{equation}\label{mc} \mathrm{Vol}(\mathcal{M}_{0,n}) = C_{n-3} \sum_{\mathcal{P}} (-1)^{|\mathcal{P}|+1} (|\mathcal{P}| -3)! \prod_{B \in \mathcal{P}} \mathrm{max}\left( 0, 1 - \sum_{i \in B} \mu_i \right)^{|B|-1},
 \end{equation}
 where $\sum_{i=1}^n \mu_i =2,$  $0<\mu_i <1, \ \forall i \in \{1,2,...,n\},$ and $C_{n} = \frac{(-4 \pi)^n}{(n+1)!}.$ Here $\mathcal{P}$ ranges over all partitions of the indices $(1, . . . ,n)$ into blocks $B.$
 \end{teorema}

From the above theorem, we can seee that the only partitions making sense are those whose size is greater or equal than three. Therefore, in this considered case, are those of size $3,$ and $4.$ Recall that the size of a partition $\mathcal{P}$ equals the number of blocks $B$ of the same partition. Therefore, if the partition has size $4$ then there is only one partition with four blocks each of which has size one. If the size of the partition is 3, then we have six different partitions of two blocks of size one and one block of size two.
From \ref{mc2}, we have $C_1 = -2 \pi,$ set $\alpha_B = \mathrm{max}\left(0,1-\sum_{i \in B}d_i \right).$ If a block has size one then $\alpha^{|B|-1}=1.$ The contribution to the sum from the partition of size 4 is just $-2\pi,$ and the contribution to the sum from the partitions of size 3 is the sum of all $d_i `$s whose indexing set is a block of size two. Namely $-2 \pi \sum_{B, \vert B \vert =2} \alpha_B.$ Then, as a final result we get

\begin{equation} \label{mc1} \mathrm{Vol}(\mathcal{M}_{0,4}) = -2 \pi \sum_{I \subset \{1,2,3,4\}: \vert I \vert = 2} (1- \sum_{i \in I} d_i)
\end{equation}

\begin{osservazione} The same result given in Equation \ref{mc1} can be achieved from a result of Thurston \cite{Th98} which shows that the general moduli space of marked curves of genus $0,$ with the additional condition that the weights sum up to $2$, is a complex hyperbolic cone manifold. Indeed, this was the key observation used by McMullen computing the volume. This complex hyperbolic manifold, in general, have a natural \emph{stratification}. The \emph{solid angle} of a cone manifold is constant on each strata. We know that for a closed complex hyperbolic manifold $M$ of dimension $n,$ its volume can be expressed by its Euler characteristic 
\begin{equation}\label{hypvol} \mathrm{Vol}(M) = \frac{(-4 \pi)^n}{(n+1)!} \chi(M)
\end{equation}

The dimension of $\mathcal{M}_{0,4}$ is $1,$ and it is homeomorphic to $\mathbb{P}^1$. Considering the analogy with the $2$-sphere with four marked points with weights $d_1,d_2,d_3,d_4,$ respectively, we can associate to these latter a cell complex, where its skeleton consist of four distinct points and six triangles. These latter correspond to the four points partitions of size two and three. With this in mind, and by applying \ref{hypvol} we get Equation \ref{mc1}.
\end{osservazione}

In the case of four points in $\mathbb{P}^1$ the result provided by Mandini in \cite{Man08} becomes

\begin{equation} \label{man4} \mathrm{Vol}(M_d) = -\pi \sum_{k=0}^{3} (-1)^{k}  \sum_{I \subset \mathcal{I}_+, |I|=k}  \left( D_I^+ - D_{I}^-\right). \end{equation} 
Now assume that $\sum_{i=1}^4 d_i =2.$ Unravelling the definition of $D_I^+$ and $D_I^-$,  we have

\begin{align*} \sum_{i \not \in I} d_i - \sum_{i \in I} d_i &= 2 -\sum_{i \in I} d_i - \sum_{i \in I} d_i \\ &= 2(1 - \sum_{i \in I} d_i). \end{align*}
By substituting in Equation \ref{man4}, we find

\begin{equation} \label{mancy} \mathrm{Vol}(M_d) = -2 \pi \sum_{k=1}^4 (-1)^{k} \sum_{I \subset\mathcal{F}, |I|=k} \mathrm{max}\left(0,1- \sum_{i \in I} d_i\right)
\end{equation}

\begin{proposizione} The volume function of the moduli space of four points in $\mathbb{P}^1$ is continuous when passing from the Fano geometry to the Calabi-Yau geometry.
\end{proposizione}

\begin{dimostrazione} All we need to show is that Equation \ref{mancy} equals Equation \ref{mc1}. When looking at Equation \ref{mancy} the subsets $I$ of size one will contribute $-2$ to the sum. Indeed, by assumption 
$d_i \in (0,1) \cap \mathbb{Q}, \forall i \in \{1,2,3,4\},$ 
therefore $\alpha_I = 1-d_i, \forall i \in \{1,2,3,4\}.$ 
By taking the sum of the four subset of size one, and observing that the sum of all the $d_i$ is by assumption two, we have 

$$ - \pi \sum_{i=1}^4 (1-d_i) =-2 \pi(4 -2) = -2 \pi. $$

The subsets of size three have zero contribute to the sum. Indeed, without loss of generality, suppose $I=\{1,2,3\},$ then $1-d_1-d_2-d_3 = 1 - 2 + d_4 = d_4 -1.$ Therefore,

$$\alpha_{\{1,2,3\}}= \mathrm{max}\left(0, d_4 -1 \right) =0.$$

Hence, it remains to sum the contributions coming from the subsets of size $2.$ But, putting all together we have

$$\mathrm{Vol}(M_d) = -2 \pi \sum_{I \subset \{1,2,3,4\}: \vert I \vert = 2} (1- \sum_{i \in I} d_i). $$

That proves the claim. \hfill $\Box$

\end{dimostrazione}

This proves, in the four point case, that the \emph{volume function} is continuous when passing from the Fano geometry to the Calabi-Yau geometry. When $\sum_{i=1}^4d_i > 2$ we may loose this behaviour, since we do not have a GIT quotient, therefore the previous techniques for calculating the volume can not be applied. The work of Alexeev \cite{Ale08} describes the moduli space of hyperplane arrangements by generalising the work of Hasset \cite{Has03} so that the one dimensional case coincides with the \emph{Hasset moduli space}. The moduli space of weighted hyperplane arrangements of dimension $1,$ is the moduli space of points in $\mathbb{P}^1.$  Since the dimension of the considered moduli space is $1,$ for any choice of weights, then we can calculate the volume in the general type geometry, by computing the degree of the log CM line bundle. To do so,  we notice that $M_d$ has a \emph{universal family} obtained by considering the product $\mathbb{P}^1 \times \mathbb{P}^1$,  fixing three points, $0,1,\infty$, and allowing the last one to move in the diagonal $\Delta = \{ (s,t) \in \mathbb{P}^1 \times \mathbb{P}^1 \ | \ s=t \ \}.$ We set $$\mathcal{X} =  \mathbb{P}^1 \times \mathbb{P}^1, \ \mathrm{and} \  \mathcal{D} = d_1 \cdot \mathbb{P}^1  \times \{0\} + d_2 \cdot \mathbb{P}^1 \times \{1\} + d_3 \cdot \mathbb{P}^1 \times \{ \infty\} + d_4 \cdot \Delta.$$ Since we fixed three points, then as divisors we can assume that they have the same hyperplane class, that will be called $h_2$. Hence, we can rewrite $\mathcal{D}$ as follows
$$ \mathcal{D} = \sum_{i=1}^4 d_i \mathrm{pr}_2^{*}h_2 + d_4  \mathrm{pr}_1^{*}h_1.$$
Because of the construction of the Hasset's moduli space we must distinguish several situations. Heuristically, we have the following phenomena: when the point on the diagonal with weight $d_4$ meets one, two or three points, then a \emph{wall} is crossed and therefore the universal family becomes singular at that point. This phenomena is encoded in the notion of \emph{stability} given in \cite{Has03} that we recall in the following

\begin{definizione} Let $\pi : (C, s_1, ..., s_n) \rightarrow B$ be a proper and flat morphism of nodal curves of arithmetic genus $g, $ and $s_1, ..., s_n$ are the sections of $\pi$ corresponding to the marked points on $C.$ A collection of data  $(g, \mathcal{A})=(g, a_1, ..., a_n)$ consist of an integer $g \geq 0$ and weights $(a_1, ..., a_n) \in \mathbb{Q}^n,$ such that $0<a_j \leq 1, \forall j \in \{1,2, ...,n\}$ and $2g - 2 + a_1 + a_2 + \cdots + a_n >0.$ We say that $\pi$ is \emph{stable} if the following conditions are satisfied
\begin{enumerate}

\item  The sections $s_1, ..., s_n$ are in the smooth locus of $\pi,$ and for every subset $\{s_{i_1}, ... s_{i_r}\}$ with nonempty intersection we have $a_{i_1} + ... + a_{i_r} \leq 1.$
\item $K_{\pi} + \sum_{i=1}^n a_i s_i$ is relatively ample.
\end{enumerate}
\end{definizione}

Clearly, in our case $g=0,$ and $M_{0,4} \simeq \mathbb{P}^1.$ The above definition, and the heuristic are rather \emph{intuitive} as the following picture shows

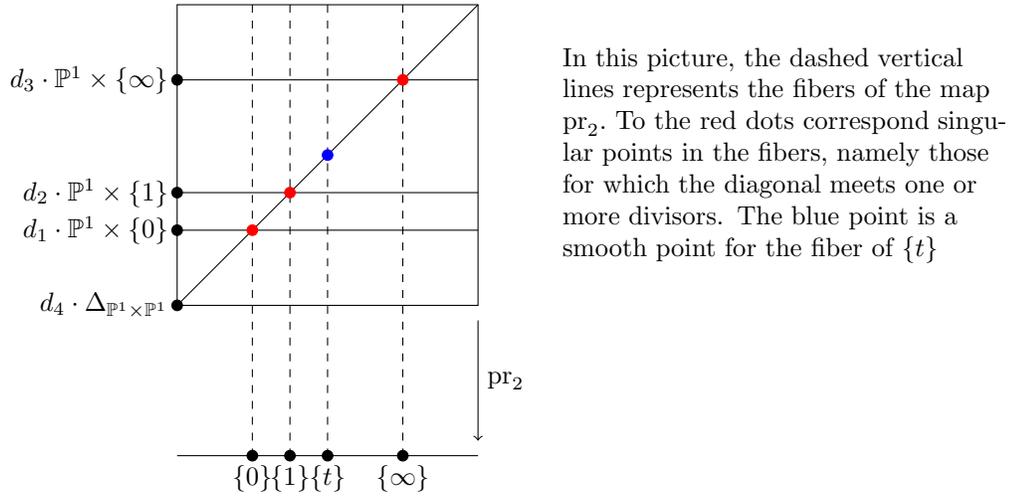
\begin{center}
\begin{tikzpicture}
\draw (0,0) -- (4,0) -- (4,4) -- (0,4) -- (0,0) ;
\draw (0,-2) -- (4,-2); 
\draw (0,0) -- (4,4);
\draw (0,1) -- (4,1);
\draw (0,1.5) -- (4,1.5);
\draw (0,3) -- (4,3);
\draw[dashed]  (1,4) -- (1,-2);
\draw[dashed] (1.5,4) --(1.5,-2);
\draw[dashed] (3,4) -- (3,-2);
\draw[dashed] (2,4) -- (2,-2);
\draw[->]  (4,-.2) -- (4,-1.8); 
\node[right] at (4,-1) { $\mathrm{pr}_2$};
\filldraw (0,0) circle (2pt) node[align=left,  left] {$d_4 \cdot \Delta_{\mathbb{P}^1 \times \mathbb{P}^1}$};
\filldraw (0,1) circle (2pt) node[align=left, left] {$d_1 \cdot \mathbb{P}^1 \times \{0\}$}; 
\filldraw (0,1.5) circle (2pt) node[align=left, left] {$d_2 \cdot \mathbb{P}^1 \times \{1\}$};
\filldraw (0,3) circle (2pt) node[align=left,left]{$d_3 \cdot \mathbb{P}^1 \times \{\infty\}$};
\filldraw (1,-2) circle (2pt) node[below] {$\{0\}$};
\filldraw (1.5, -2) circle (2pt) node[below] {$\{1\}$};
\filldraw (2, -2) circle (2pt) 
node[below] {$\{t\}$};
\filldraw (3,-2) circle (2pt) node[ below] {$\{\infty\}$};
\filldraw[red] (1,1)  circle (2pt);
\filldraw[red] (1.5, 1.5) circle (2pt);
\filldraw[red] (3,3) circle (2pt);
\filldraw[blue] (2,2) circle (2pt);
\node[text width=6cm, anchor=west, right] at (5,2)
    {In this picture, the dashed vertical lines represents the fibers of the map $\mathrm{pr}_2.$ To the red dots correspond singular points in the fibers, namely those for which the diagonal meets one or more divisors. The blue point is a smooth point for the fiber of $\{t\}$ };
\end{tikzpicture}
\end{center}

When looking at the above picture, we shall assume that the weight $d_4$ is the smallest among the weights. So, while $d_4$ is free to move, and meet one (or more) points in the diagonal, then the rest of the points remain fixed. 

\begin{itemize}

\item  The following inequalities holds, when the diagonal does not meet any of the fixed points \begin{equation} \label{casouno}
 \begin{cases}
 d_i + d_j >1 & \ \text{$i,j \in \{1,2,3\}$} \\
 d_k + d_4 < 1 & \ \text{for all $k \in \{1,2,3,4\}$}
 \end{cases}
 \end{equation}\end{itemize}

\begin{itemize}

\item  If $d_4$ grows, then it meets one (or more) points along the diagonal, and a desingularization is needed  \begin{equation}
\label{casodue}
 \begin{cases}
 d_i + d_4 >1 & \ \text{$ i \in \{1,2,3\}$} \\
 d_j + d_4 \leq 1 & \ \text{$i \neq j$.}
 \end{cases}
 \end{equation}\end{itemize}
 
As the above picture suggest we can take the projection onto the second factor of $\mathcal{X}$, namely $\mathrm{pr}_2 : (\mathcal{X}, \mathcal{D}) \rightarrow \mathbb{P}^1,$ to get a  $\mathbb{Q}$-Gorenstein family of log general type varieties. The stability conditions of Hasset suggest that we should choose, as a relatively ample line bundle for $\mathrm{pr}_2,$ the log canonical polarization $K_{\mathcal{X}/\mathbb{P}^1} + \mathcal{D}.$

We begin by studying the \ref{casouno}. We chose $K_{\mathcal{X}/\mathbb{P}^1} + \mathcal{D}$ as a relative polarization for the family $\mathrm{pr}_2 \colon (\mathcal{X}, \mathcal{D}) \rightarrow \mathbb{P}^1.$ By a direct application of Corollary \ref{loggtcr} we have that 

\begin{equation} \label{cmgtliscio} c_1 (\lambda_{\mathrm{CM}, \mathcal{D}})= 2 d_4 \left( \sum_{i=1}^4 d_i -2 \right). \end{equation}In \ref{casodue}, it is the situation when the diagonal meets one point. Suppose, without loss of generality, that the diagonal meets the zeroth fiber. The universal family becomes singular at that point, therefore we take the blowup at $\{0\}$ of the total space of the universal family. Set $\tilde{\mathcal{X}}=\mathrm{Bl}_0 (\mathcal{X}).$ The divisor $\mathcal{D},$ modifies as follows

\begin{align*} \mathcal{D} &= d_1 (\mathrm{pr}_2^{*}h_2 -E) + (d_2 + d_3) \mathrm{pr}_2^{*}h_2 + d_4 (\mathrm{pr}_2^{*}h_2 + \mathrm{pr}_1^{*}h_1 -E) \\ &= \sum_{i=1}^4 d_i \mathrm{pr}_2^{*}h_2 + d_4 \mathrm{pr}_1^{*}h_1 - (d_1 + d_4) E.
\end{align*}Then, 

$$ K_{\tilde{\mathcal{X}}/\mathbb{P}^1} + \mathcal{D} = \left(\sum_{i=1}^4 d_i -2 \right) \mathrm{pr}_2^{*}h_2 + d_4 \mathrm{pr}_1^{*}h_1 - (d_1 + d_4 +1)E. $$
By applying directly Corollary \ref{loggtcr} we find Equation \ref{cmgtliscio}. Indeed the cohomology ring of the blow up suggests that the intersection products $\mathrm{pr}_i^{*}h_i \cdot E =0, i=1,2$. Moreover, the pushforward of the constant term $(d_1 + d_4 +1)$ coming from $E^2=-1$ term, is zero, therefore the only term that survives is only $\mathrm{pr}_{2*} 2d_4 \left( \sum_{i=1}^4 d_i -2 \right) \mathrm{pr}_2^{*}h_2 \cdot \mathrm{pr}_1^{*}h_1$, that yields to \ref{cmgtliscio}. The same holds if more than one desingularization is needed. The formula of the volume function for the log general type case is therefore constant with respect to the log canonical polarization. We notice also that when Equation \ref{cmgtliscio} is evaluated in the Calabi-Yau zone then the volume function is zero. With respect to the anticanonical log polarization given in the Fano case we have a change of sign. This shows that the volume function has a discontinuity point in the Calabi-Yau zone. However, according to \cite[Theorem 5.5.2]{Ale15}, and \cite[Theorem 1.5]{Ale08}, if we choose weights whose sum is slightly greater than $2$ the moduli space do not change, and it is described as the GIT quotient \ref{git4}. Then, at least for \ref{casouno}, we can chose as relative polarization $-K_{\mathcal{X}/\mathbb{P}^1}$. Therefore, a fast computation proves that 

$$ c_1 (\lambda_{\mathrm{CM}, \mathcal{D}})= 4 d_4. $$ We can easily observe that up to a normalization factor, $\pi$, it coincides with the result obtained by applying \ref{mc1}, and \ref{mancy} with the following order of weights

$$ d_i < \sum_{j \neq i} d_j, \ \forall i \in \{1,2,3,4\}, $$

\begin{align}  d_4 + d_2 < d_3 + d_1 \nonumber \\
d_4 + d_3 < d_2 + d_1  \nonumber \\ 
d_4 + d_4  < d_2 + d_1. \nonumber
\end{align}

We resume all these results in the following

\begin{proposizione} \label{propordine} Let $M_d$ be the moduli space of weighted hyperplane arrangements of dimension $1.$ Suppose that the weights $d_i \in (0,1) \cap \mathbb{Q}, i \in \{1,2,3,4\}$ satisfy the following conditions
\begin{itemize}

\item $d_j < \sum_{i \neq j} d_i, \forall j \in \{1,2,3,4\};$ 
\item  \begin{align}  d_4 + d_2 < d_3 + d_1 \nonumber \\
d_4 + d_3 < d_2 + d_1  \nonumber \\ 
d_4 + d_4  < d_2 + d_1. \nonumber
\end{align} \end{itemize}

Then, for small weights the volume of $M_d$ changes continuously along the Fano, Calabi-Yau and general type geometry and it is given by $ 4 \pi d_4.$

\end{proposizione} 

\begin{osservazione} In the condition of Proposition \ref{propordine}, when the sum of the weights largely exceed two, then we must chose another polarization, i.e. the log canonical polarization, as the Hasset compactification suggests. The calculations shows that the volume of $M_d$ in this case is given by $2 d_4 \left( \sum_{i=1}^4 d_i -2 \right).$ Note that when approaching the Calabi-Yau geometry from the \emph{far away} general type geometry, then the volume 
of $M_d$ goes to zero.
\end{osservazione}

 \section{The general case}
 
 In order to prove Theorem \ref{importante}, we recall the following definition
 
 \begin{definizione}\label{logwp} \cite{Tam20} Let $f: (\mathcal{X}, \mathcal{D}) \longrightarrow B$ be a family of log K-polystable Fano varieties of pure dimension $n.$ Assume that $\mathcal{D}:=\sum_k (1-\beta_k)[s_k=0]$ is simple normal crossing. Then, we define the log Weil-Petersson metric as the following fiber integral
 \begin{equation} \label{eq1} \omega_{\mathrm{WP}}^0 := - \int_{\mathcal{X}^0 / B^0} \omega_{\mathcal{X}^0}^{n+1}.
 \end{equation}
\end{definizione}

In \cite[Theorem 2.3]{Tam20} it is proven that Equation \ref{eq1} can be extended to the whole family $f: (\mathcal{X}, \mathcal{D}) \longrightarrow B,$ and it extends as a curvature for a metric on the descended log CM line bundle $\lambda_{\mathrm{CM},\mathcal{D}}.$

\begin{dimostrazione}(Of Theorem \ref{importante}). Consider $n$ fixed points $\{p_i\}_{i=1}^n$ in general position in $\mathbb{P}^1.$ Each of these points has weight
$d_i \in (0,1)\cap \mathbb{Q},$ where $d_i = 1-\beta_i$ is related to the conic angle $\beta_i.$ Fix a conic KE metric in $\mathbb{P}^1,$ $\omega_d \in c_1 (\mathcal{O}(1)),$ where $d=(d_1, ...,d_n).$ Note that when $d=(1, ..., 1)$ then $\omega_d$ coincides with the Fubini-Study metric on $\mathbb{P}^1.$ When $\sum_{i=1}^n d_i =2$ then $\omega_d$ is a Calabi-Yau metric, that is

\begin{equation}\label{calabi} 
\omega_d = i \Omega \wedge \overline{\Omega}
    \end{equation}
    
Where, by choosing an affine coordinate $z$ for the points $\{p_i\}_{i=1}^n$, such that $z(p_i)=c_i \in \mathbb{C}, \forall i \in \{1,2, ...,n\},$ the holomorphic $1$-form $\Omega$ is given by 

$$ \Omega = \frac{dz}{\prod_{i=1}^n (z - c_i)^{d_i}}. $$

In \cite[section 4.1]{dBSp18a} it is proven that \ref{calabi} is indeed the metric used by McMullen in \cite{MM00}, that is the Weil-Peterson metric. Namely,

$$ \omega_{WP} = - i \partial \overline{\partial} \log \int_{\mathbb{P}^1} i \Omega \wedge \overline{\Omega}. $$

Now, let $f: (\mathcal{X}, \mathcal{D}) \longrightarrow \mathbb{P}^1$ be a familily of $K$-polystable log Fano hyperplane arrangements of dimension one, like in section 2.2.  That is, $\mathcal{X}= \mathbb{P}^1 \times \mathbb{P}^1,$ $\mathcal{D}=\sum_{i=1}^{n-1} d_i \mathrm{pr}_2^{*} h + d_n (\mathrm{pr}_1^{*}l + \mathrm{pr}_2^{*}h).$ Then $f$ is the projection map onto the second factor of $\mathcal{X}.$ Furthermore, we assume for this family that the sum of all the weights is two. We choose as a relative ample line bundle for $f$ $\mathcal{L}= - K_{\mathcal{X}/\mathbb{P}^1}.$ Fiberwise we have 
$$ - i \partial \overline{\partial} \log \int_{\mathcal{X}_t} i \Omega_t \wedge \overline{\Omega_t} = - i \partial \overline{\partial} \log \int_{\mathcal{X}_t}  \omega_{d,t}^{n+1}. $$

Then, because of Definition \ref{logwp} the above is exactly the fiberwise definition of the log Weil Petersson metric, and hence the curvature for a metric on $\lambda_{\mathrm{CM},\mathcal{D}}.$ That means,

\begin{equation}\label{equal} \mathrm{Vol}(\mathcal{M}_{0,n}) = \mathrm{Vol}(\lambda_{\mathrm{CM}, \mathcal{D}}, \omega_{WP}). \end{equation} 

Because of the second assertion of Proposition \ref{prop3}, we see that the log CM line bundle it is a multiple of the anti-canonical bundle of $\mathbb{P}^1.$ Therefore, the volume can be achieved with the Jeffrey-Kirwan residue theorem, that leads to equation in \cite{MM00} (Theorem 8.1) as wanted. \hfill $\Box$
\end{dimostrazione}

\begin{osservazione} A combinatorial proof like in the case of four points it is rather complicated. As a check we provide a Python script that shows numerically the equality of the two formulas, that can be found in APPENDIX A.

\end{osservazione}

 \section{Final remarks}
 
 We wish to conclude this section with an observation. Consider the log Fano pair $(\mathbb{P}^1, \sum_{i=1}^4 d_i p_i)$.  Chose some coordinate $z \in \mathbb{P}^1,$ such that for a fixed $i \in \{1,2,3,4\}$ we have $z(p_i)=0$. Consider the vector field $v = z \partial_z,$ and notice that it generates a one parameter subgroup $\lambda: \mathbb{C}^{*} \hookrightarrow \mathrm{GL}_1 (\mathbb{C})$, that acts on $\mathbb{P}^1$ in the following fashion $(\lambda(t), z):= t \cdot z.$  
This actions easily translates on the divisor, $\mathcal{D}_t := t \cdot (\sum_{i=1}^4 d_i p_i)$, and on the anticanonical polarization of $\mathbb{P}^1,$ namely $\mathcal{O}_{\mathbb{P}^1}(2)$. This data define a test configuration $\{(\mathbb{P}^1_t , \mathcal{D}_t, -K_{\mathbb{P}^1_t} )\}_{t \in \mathbb{C}^{*}}=(\mathcal{X}, \mathcal{D}, K_{\mathcal{X}})$ whose central fiber $\mathcal{X}_0$ is just a point and the divisor, at the central fiber, behave as

\begin{equation} \label{centralfiber} \mathcal{D}_0 = d_i \{\infty\} + \sum_{i \neq j} d_j \{0\}. \end{equation}
Then, this latter defines an \emph{integral test configuration}. As we mentioned in previous discussions and in \cite{Fuji19a} there exist a bijection between integral test configurations and \emph{dreamy prime divisors}. We observe that since we are dealing with toric varieties, then every prime divisor is also dreamy (\cite{Fuji19a}). In particular in our central fiber \ref{centralfiber} we have the sum of two dreamy prime divisor 

$$ W_1 = \{0\}, \ \mathrm{and} \ W_2= \{\infty\}. $$We have that
$$\mathrm{ord}_{W_1} ( \mathcal{D}_0)= d_i, \ \mathrm{and} \ \mathrm{ord}_{W_2}( \mathcal{D}_0) = \sum_{i \neq j} d_j$$ 
The divisor $W_1$ leads to $\mathrm{DF}_{\mathcal{D}} (\mathcal{X}, \mathcal{L}) = \sum_{i \neq j}d_j - d_i.$ We will now show that the volume of $M_{0,4}$ can be obtained by summing all the Donaldson-Futaki invariant of the corresponding integral test configuration. In \ref{man4} every single term $D_F^+ - D_F^-$ , by the Jeffrey-Kirwan \emph{residue theorem} \cite{JK95} is associated to 

\begin{equation} \label{pezzo} \frac{i_B ^{*} c_1 (\mathcal{L}_d)(X)}{e_B (X)} = (-1)^{4-k} \mu(B) X^{-3}
\end{equation}
where $B \in \mathcal{F} = \{ (z^1, ..., z^4) \in (\mathbb{P}^1)^4 | \ z^j \in \{0, \infty \} \}$ is a fixed point for the maximal torus action. In \ref{pezzo} we used the fact that since $\mathcal{L}_d$ is a prequantum line bundle and the maximal torus is the unit circle then by \cite{Jef19}, Lemma 9.31, the image of the moment map $\mu(B)$ is given by the weighted sum of the height function on each $\mathbb{P}^1,$ namely

$$ \mu(B) = \sum_{i=1}^4 d_i \mu^i(z^i), $$

where $\mu^i (z^i)$ is the height function, i.e. the moment map for the maximal torus action on each $\mathbb{P}^1.$ Clearly, the fixed points for this action are the standard basis element $e_1, e_2$ of $\mathbb{C}^2$ that corresponds to $0, \infty$ respectively. We have 

\begin{displaymath}
\mu^i (e_i)=
 \begin{cases}
 1  &  \ \text{if $i=1$} \\
 -1 & \ \text{if $i=2.$}
 \end{cases}
 \end{displaymath}
 
 It follows in \ref{man4} that the index $F$ is the one that tells how many $e_2$ we have in the fixed point $B.$ 
  
  If $|F| \geq 2$ we can always think of a integral test configuration for which two points came together into one point, where the weight of the latter is the sum of two or more weights. With this in mind, it is easy to convince ourselves that the size of the index set $F$ tells how many points should come together into one point, and hence how many weights should be summed at that point. Therefore, the general term 
  
  $$ \frac{i_B^{*} c_1(\mathcal{L})}{e_B (X)} = - \frac{(-1)^{4-k}}{x^3} \mathrm{DF}_{\mathcal{D}} (\mathcal{X}, \mathcal{L})(B),$$
  as wanted. This observation, can be immediately generalised for a configurations of $n-$points in $\mathbb{P}^1.$ The key point is that in the work of \cite{Fuji19b} it is shown that in the case of log Fano hyperplane arrangements the notion of K-stability coincides with the notion of GIT-stability. Moreover by changing the weights, the corresponding moduli spaces are \emph{well behaved} so that the Jeffrey-Kirwan localization formula of \cite{JK95} can be applied for calculating the volume. The GIT quotient describing the hyperplane arrangements of dimension one is the following
  
  $$ M_d = (\mathbb{P}^1)^n //_{\mathcal{L}_d} \mathrm{SL}_2 \mathbb{C}. $$
  The dimension of $M_d$ is $n-3.$ We apply directly the Jeffrey-Kirwan theorem for the case of $\mathrm{SU}(2),$ on the top class namely
  $$ k(c_1(\mathcal{L}_d)^{n-3}) e^{\omega_0}[M_d ] = \frac{n_0}{2} \mathrm{Res}_{X=0}\left( 4X^2 \sum_{F \in \mathcal{F}_+} \int_{F} \frac{i_F ^{*} c_1(\mathcal{L}_d)^{n-3}(X) } {e_F (X)}dX \right) $$
  
   since $\mathcal{L}_d$ is a prequantum line bundle and the maximal torus is the unit circle then by Lemma 9.31 of \cite{Jef19} $i_f^{*} c_1(\mathcal{L}_d)^{n-3}(X) = (-1)^{n-3} \mu(F)^{n-3}(X).$ Hence,

   $$ \frac{i_F ^{*} c_1 (\mathcal{L}_d)^{n-3}}{e_F (X)} = -(-1)^{n-k} (\mu(F)(X))^{n-3} X^{-n-3}. $$
   Because of the above observation, $\mu(F) = \mathrm{DF}_{\mathcal{D}}(\mathcal{X}_F, \mathcal{L}_F)= D_F^{+} - D_F^{-}$.  We have 
   
   \begin{equation} \label{pezzobis} 4X^2 \frac{i_F ^{*} c_1 (\mathcal{L}_d)^{n-3}}{e_F (X)} = -(-1)^{n-k} \mathrm{DF}_{\mathcal{D}}(\mathcal{X}_F, \mathcal{L}_F)^{n-3} X^{-1}. \end{equation}
   The residue of $X^{-1}$ in \ref{pezzobis} is one. By taking the sum, we proved the following result
  
  \begin{teorema} \label{risfut} The volume of the moduli space of log Fano hyperplane arrangements of dimension one is the sum of the Donaldson-Futaki invariants associated to test configurations $(\mathcal{X}_F , \mathcal{L}_F).$ Where $F \in \mathcal{F}_+$ is like in \cite{Man08}, and $(\mathcal{X}_F, \mathcal{L}_F)$ is the test configuration for which $|F|$ points come together into one point.
  
  $$ \mathrm{Vol}(M_d) = -\frac{(2 \pi)^{n-3}}{2(n-3)!} \sum_{i=0}^{n-3} (-1)^k \sum_{F \in \mathcal{F}_+, \ |F|=k} \mathrm{DF}_{\mathcal{D}}(\mathcal{X}_F, \mathcal{L}_F)^{n-3}. $$
 \end{teorema}
 
 \begin{osservazione} It is natural to ask if more generally one can express the volumes as \emph{sums} of CM degrees of special unstable families, using the $\Theta$-stratification \cite{HL} to generalise \emph{Kirwan's} conditions \cite[Theorem 5.6]{HL}. 
 \end{osservazione}

\section*{APPENDIX A}

A Python script for calculating the equality of the two mentioned volume formulas is given below.

\lstset{language=Python}
\lstset{frame=lines}
\lstset{caption={Test of Equivalence}}
\lstset{label={lst:code_direct}}
\lstset{basicstyle=\footnotesize}
\begin{lstlisting}


import numpy as np
import math
from fractions import Fraction
import sympy
import itertools
from itertools import combinations
import random

def subsets_k(collection, k): 
    yield from partition_k(collection, k, k)
    
def partition_k(collection, min, k):
  if len(collection) == 1:
    yield [ collection ]
    return

  first = collection[0]
  for smaller in partition_k(collection[1:], min - 1, k):
    if len(smaller) > k: continue
    # insert `first` in each of the subpartition's subsets
    if len(smaller) >= min:
      for n, subset in enumerate(smaller):
        yield smaller[:n] + [[ first ] + subset]  + smaller[n+1:]
    # put `first` in its own subset 
    if len(smaller) < k: yield [ [ first ] ] + smaller


def mandini_volume(V):
    
    '''Computes the Volume of the Moduli space of poligons
    
       :params V, list of integers'''
    
    
    
    N = len(V)
       
    num = Fraction(np.power(2, 2*N-7))
    Prefactor = Fraction(num, math.factorial(N-3))*np.power(sympy.pi, N-3)
    
    Result = 0
    
    for I in range(0,N):
        R = list(combinations(V,I))
        factor = 0
        for j in range(0, np.array(R, dtype = 'object').shape[0]):
            Q = 1
            S = np.power(np.max([0, 1- np.sum(R[j])]), N-3)
            Q *=S
            factor += (-1)**(I+1)*Q
        Result += Fraction(factor)
    return (Result*Prefactor)

def mcmullen_volume(V):
    
    '''Computes the volume of the moduli space of gurves of genus 0
       
       :params V, list of integers'''
    
    N = len(V)
    
    Prefactor = Fraction((-4)**(N-3),math.factorial(N-2)) * sympy.pi**(N-3)

    Result = 0
    
    for P in range(3,N+1):
        R = list(subsets_k(V, P))
        factor = 0    
        for j in range(0,np.array(R, dtype="object").shape[0]):
            Q = 1
            for i in range(0,P):
                B = len(R[j][i])
                S = np.max([0,1-np.sum(R[j][i])])**(B-1)
                Q *=S
            factor += (-1)**(P+1) *math.factorial(P-3)*Q
        Result += factor
        
    return (Result*Prefactor)

class AnomalyTest():
    
    '''
        Anomaly test, is a class that performs a test to check wether the functions mandini_volume and mcmullen_volume achieve the same results
       
       :params length, is an integer representing the number of marked points in the moduli space
       :params total_tests, is an integer representing the number of tests to be achieved
       
       :method random_vector, returns a test vectors of integers
       :method run_test, performs the test. 
       
    '''
    
    def __init__(self, length, total_tests):
        
    
        self.length = length    
        self.total_tests = total_tests
        
        if not isinstance(self.length, int):
            raise ValueError(f'lenght must be integer, found {self.length}')
            
        elif not self.length >= 3:
            raise ValueError(f'length must be greater or equal than 3, found {self.length}')
        
        elif not isinstance(self.total_tests, int):
            raise ValueError(f'total_test must be integer, found {self.total_test}')
            
        elif self.total_test < 0:
            raise ValueError(f'total_test must be positive, found {self.total_test}')
        
        
        self.random_vector()
        
        
        
    def random_vector(self):
        
        V = []                                 
    
        for i in range(0,self.length):
            n = random.randint(1,30)
            V.append(n)
        
        V = list(np.asarray(V)*Fraction(2)/np.sum(V))
        
        self.test_vector = V
    
        
        
    def run_test(self):
        
        
        test_iterations, anomalies = np.zeros((1,), dtype=int), np.zeros((1,), dtype=int)

        while test_iterations < self.total_tests: 
            
            V = self.test_vector

            check= any(entry > 1 for entry in V)

            if check: 
                continue

            test_iterations+=1

            Vol_1 = mandini_volume(V)
            Vol_2 = mcmullen_volume(V) 
            
            if Vol_1 == 0 and Vol_2 == 0: 
                Ratio = 1.0
            else:
                Ratio = float(Vol_1 / Vol_2) 

            if Ratio != 1: 
                anomalies += 1
                print (V)
                print('Vol_1, Vol_2:', Vol_1, Vol_2)
                print('The ratio of the two values of Vol(M_d):', r) 

        print('Total anomalies', anomalies)   




\end{lstlisting}



\end{document}